\documentclass[12pt]{article}

\usepackage{epsfig}

\makeatletter
\renewcommand{\@biblabel}[1]{#1.}
\makeatother \textwidth160mm \textheight250mm
\newcommand{\RR}{{\rm I\kern-0.14em R }}
\newcommand{\NN}{{\rm l\kern-0.14em N}}
\newcommand{\CC}{{\rm\raise 0.192ex\vbox{\hrule height
                  1.22ex width 0.8pt}\kern-0.29em C}}

\begin{document}

\begin{center} \bf \Large
About the group law for the Jacobi variety
of a hyperelliptic curve
\end{center}


\begin{center}
{\bf Frank Leitenberger } \\
{\small \it Fachbereich
Mathematik, Universit\"at Rostock, Rostock, D-18051, Germany. \\
e-mail: frank.leitenberger@mathematik.uni-rostock.de }
\end{center}


\small \normalsize

\vspace{0.2cm}

\begin{quote}\small {\bf Abstract.}
We generalize the group law
of curves of degree three by chords and tangents
to the Jacobi variety of a hyperelliptic curve.
In the case of genus 2
we accomplish the construction by a cubic parabola.
We derive explicit rational formulas for the addition
on a dense set in the Jacobian.
\end{quote}

\vspace{0.3cm}

\noindent
{\bf 1. Introduction}

\vspace{0.1cm}

\noindent
The intention of this remark is an explicit description
of the group law of hyperelliptic curves.
It appears that it is possible to generalize
the chord and tangent method
for curves of degree three in a very naive way
by replacing points by point groups of $g$ points
and by replacing lines by certain interpolation functions.

Explicit descriptions of the group law play a less important role
in the history of the subject. They appear first in the new literature.
Cassels remarked 1983: "{\it I cannot even find in the literature
an explicit set of equations  for the Jacobian of a curve of genus
{\rm 2} together with explicit expressions for the group operation
in a form  amenable to calculation ...}" (cf. \cite{Cas,Gr}).
Mazur remarked 1986: "{\it ... a naive
attempt to generalize this group structure {\rm [of degree 3
curves]} to curves of higher degree (even quartics) will not
work.}" (cf. \cite{Ma}, p. 230).
With the development of cryptography arose algorithms for the group
law. In 1987 Cantor described the group law of a hyperelliptic
curve in the context of cryptography (cf. \cite{Ca,K}). Later
group laws of more general classes of curves were described in
\cite{Hu,Vol}. These group laws work step by step and do not
allow a visualization.

In this remark we derive explicit formulas for the group law
for the Jacobi variety of a curve of genus 2
starting from an interpolating cubic parabola.
As the above algorithms perform the reduction in several steps
we execute the reduction in only one step.
The case $g>2$ can be performed by rational interpolation functions
analogously.
These interpolation functions were first considered by
Jacobi in connection with Abel's theorem
(cf. \cite{J}).
Our formulas are much simpler than analogous
formulas derived by Theta functions in \cite{Gr} p.114-116,
\cite{M,Y}.
A different geometric interpretation
was given by Otto Staude in \cite{St}.

\vspace{0.5cm}

\noindent
{\bf 2. Preliminaries}

\vspace{0.1cm}

\noindent
Consider a hyperelliptic curve
$C=\{\ (x,y)\in \CC^2 \ |\ \ y^2=p(x)\  \} \cup \{\infty\}$
of genus $g$
where $p(x)=a_0 x^{2g+1}+a_1 x^{2g}+\ ...\ +a_{2g+1}$ is a complex
polynomial with $a_0\neq 0$, $g\geq 1$ without double zeros.
$C$ is endowed with the involution
$\overline{(x,y)}:=(x,-y)$,
$\overline{\infty}:=\infty $.
The Jacobi variety of $C$ is the Abelian group
\[ Jac(C)=Div^0(C)/Div^P(C),\]
where $Div^0(C)$ denotes the group of divisors of degree $0$ and
$Div^P(C)$ is the subgroup of principal divisors (i.e. the zeros
and poles of analytic functions), cf. \cite{Mu}.
We find in every divisor class of
Jac(C) an unique so called reduced divisor of the form
\[n_1P_1+...+n_mP_m-(n_1+...+n_m)\infty,\]
where $n_1+...+n_m \leq g$,
$P_i\neq P_j,\overline{P_j},\infty$ for $i\neq j$
and $n_i = 1$ if $P_i = \overline{P_i}$
(cf. \cite{Mu}).
We remark that $-(P-\infty)\sim\overline{P}-\infty $ $(*)$
and $P_1+\cdots +P_h\sim h\infty$ $(**)$ if $P_1,\cdots ,P_h$
are the finite intersections of $C$ with an algebraic curve.

Now we consider the two reduced divisors
\[ J_1=P_1+ \cdots +P_{h_1}-h_1\infty ,\ \ \ \ \ \ \
   J_2=Q_1+ \cdots +Q_{h_2}-h_2\infty  \]
with $0\leq h_1,h_2\leq g$
(in this notation points $P_i, Q_j$ can occur repeatedly).
Without restriction of generality we
have $r$ $(0\leq r \leq h_1,h_2)$ pairs
$P_{h_1-k}=\overline{Q_{h_2-k}}$, $k=0,\cdots ,r-1$. Because of
$P+\overline{P}\sim 2\infty$ it follows
\[  J_1+J_2\sim P_1+\cdots +P_{h_1-r}+Q_1+
  \cdots +Q_{h_2-r}-(h_1+h_2-2r)\infty. \]
In the case $h_1+h_2-2r\leq g$ we have already
an reduced divisor on the left side.
Otherwise we consider the interpolation function
\[      y=\frac{b_0x^p +  \cdots +b_p}{
           c_0 x^q+c_1x^{q-1}\cdots +c_q}=:\frac{b(x)}{c(x)} \]
(cf. \cite{J}) with
      $p=\frac{h_1+h_2+g-2r-\varepsilon }{2}$,
      $q=\frac{h_1+h_2-g-2r-2+\varepsilon }{2}$
where $\varepsilon$ is the parity of $h_1+h_2+g$.
We have $p+q+1=h_1+h_2-2r$ degrees of freedom.
We can determine the coefficients unique up to
a constant factor so that we interpolate
the points $P_i$, $Q_j$
(in the case of a multiple point $P$ we require
a corresponding degree of contact with $C$).
These $h_1+h_2-2r$ points lie on the
algebraic curve $y c(x)-b(x) = 0$.
It follows $p(x)c^2(x)- b^2(x)=0$.
On the left side we have a polynomial of degree
$\leq h_1+h_2-2r+g$.
Therefore we obtain $h_3 \leq g$
further finite intersections $R_1,\cdots ,R_{h_3}$.
With $(*),(**)$ it follows that
\[ \overline{R_1}+\cdots + \overline{R_{h_3}}-h_3\infty   \]
is the reduced divisor for $J_1+J_2$.

It appears that only for $g=1,2$
nonfractional interpolation functions
are sufficient. Consider the case $g=2$.
Let $J_1=P_1+P_2-2\infty$, $J_2=Q_1+Q_2-2\infty$ be
two reduced divisors with $P_i\neq \overline{Q_j}$.
The interpolation polynomial
\[ y=b_0x^3+b_1x^2+b_2x+b_3 \]
through the $P_i, Q_i$ (possibly with multiplicities)
intersects $C$ for $b_0\neq 0$ in two further finite points
$R_1$ and $R_2$
with $R_1\neq \overline{R_2}$.
The result is
\[J_1+J_2=\overline{R_1}+\overline{R_2}-2\infty . \]

\vspace{0.5cm}

\epsfig{file=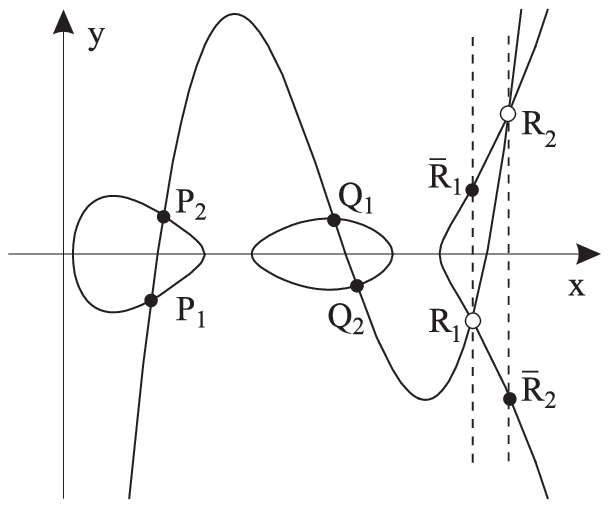, width=7.5cm}

\vspace{0.5cm}
{\small
\[  {\rm Figure\ 1.}\ \ \
(P_1+P_2-2\infty)\  +\  (Q_1+Q_2-2\infty)\ \sim\
\overline{R_1}+\overline{R_2} -2\infty     \]}

\noindent
{\it Remark:} In the real case,
contrarily to the case $g=1$ for $g=2$
the reduction of the sum of two divisors with real points
can give a sum of two complex conjugated points.

\vspace{0.5cm}

\noindent
{\bf 3. Explicit formulas}

\vspace{0.1cm}

\noindent
We use the construction in order to derive explicit formulas
in the case $g=2$.
We consider only the generic case where $b_0 \neq 0$
and all $P_1,P_2,Q_1,Q_2$
have different nonvanishing $x$-coordinates.
In this case we have the interpolation polynomial
\[  y  = b(x) = b_0 x^3+b_1x^2+b_2x+b_3
 = \sum_{i=1}^4 y_i \prod_{j\neq i} \frac{(x-x_j)}{(x_i-x_j)}.\]
For the $x$-coordinates of the intersections with the curve
$y^2= a_0 x^{5} + a_1 x^{4} +\ ...\ +a_5$
we obtain
\[  (b_0 x^3+b_1x^2+b_2x+b_3)^2-a_0 x^{5}-a_1 x^{4} -\ ...\ -a_5 =0. \]
For the six intersections it follows
\[  x_1+x_2+x_3+x_4+x_5+x_6=\frac{a_0-2b_0 b_1}{b_0^2},
\ \ \ \ \ \ \
    x_1x_2x_3x_4x_5x_6     =\frac{b_3^2-a_5}{b_0^2}.       \]
According to Vieta $x_5$ and $x_6$ are solutions of the quadratic
equation
\begin{eqnarray} \label{equation7}
   x^2+(x_1+x_2+x_3+x_4-\frac{a_0-2b_0 b_1}{b_0^2} )x
 +\frac{b_3^2-a_5}{b_0^2x_1x_2x_3x_4}=0.
\end{eqnarray}
Therefore we obtain
\[  \overline{R_1}=(x_5,-b_0 x_5^3 -b_1x_5^2-b_2x_5 - b_3),
\ \ \ \ \ \ \
\overline{R_2}=(x_6,-b_0 x_6^3 -b_1x_6^2-b_2x_6 - b_3). \]

%
%

\vspace{0.5cm}

\noindent
{\bf 4. Rational formulas}

\vspace{0.1cm}

\noindent
The group law of the previous section
contains a root operation.
It is possible to avoid roots
by the representation of divisors by Mumford and Cantor
(cf. \cite{Ca,Mu}).
We present a reduced divisor $P_1+P_2=(x_1,y_1)+(x_2,y_2)$
by the pair of polynomials
\[ ((x-x_1)(x-x_2),\frac{y_2-y_1}{x_2-x_1}(x-x_1)+y_1)=:(A(x),B(x))
= (x^2+\alpha x+\beta ,\gamma x+\delta) \] if $x_1\neq x_2$. A
divisor $2P_1=2(x_1,y_1)$ has the representation
$ ((x-x_1)^2,
\frac{{p\; }'(x_1)}{2y_1}(x-x_1)+y_1).$
The divisors of the form $D=P_1=(x_1,y_1)$ form the so called
Theta divisor $\Theta$. We can represent $(x_1,y_1)$ by
the pair  $(x-x_1,y_1)$. Now we consider
the sum
\[ (A_1(x),B_1(x)) + (A_2(x),B_2(x)) = (A_3(x),B_3(x)). \]
The coordinates
$\alpha ,\beta ,\gamma ,\delta $
form a coordinate system on $Jac(C)- \Theta$.
We show that the group law has a rational form in the generic case
$Nb_0\beta_1\beta_2 \neq 0$
(cf. below for $b_0$, $N$).


We can replace the $x_i,y_i$
of the cubic interpolation polynomial through the
$\alpha_i,\beta_i,\gamma_i,\delta_i$
by a Groebner basis calculation.
We insert the expressions for $y_i$ into $b(x)$ and we
consider the ring $\CC[x,y,a_1,a_2,b_1,b_2][x_1,x_2,x_3,x_4]$,
the order $x_1<x_2<x_3<x_4$ and the ideal
\[ (\
(x_1-x_2)(x_1-x_3)(x_1-x_4)(x_2-x_3)(x_2-x_4)(x_3-x_4)(y-b(x)),  \]
\[ \alpha_1+x_1+x_2,\alpha_2+x_3+x_4,\beta_1-x_1x_2,\beta_2-x_3x_4 \  ). \]
By a computer calculation we find the first Groebner basis element
\[ (a_1^2 - 4 b_1)(a_2^2 - 4 b_2)\
(\ ((\beta_1 -\beta_2 )^2 +(\alpha_1-\alpha_2)
    (\alpha_1 \beta_2 -\alpha_2 \beta_1)) y - \tilde{b}(x)\ ) \]
where $\tilde{b}(x)$ is independent from the $x_i$.
We require that the discriminants of $A_1,A_2$ do not vanish.
Furthermore we have
\[ b_0=\frac{1}{N}
 ((\beta_2-\beta_1) (\gamma_1 - \gamma_2) +
(\alpha_1 - \alpha_2) (\delta_1 - \delta_2)) ,\]
\[b_1=\frac{1}{N}
 ((\alpha_2\beta_2-\alpha_1 \beta_1) (\gamma_1 - \gamma_2)
  + (\alpha_1^2-\alpha_2^2-\beta_1+\beta_2 )(\delta_1 - \delta_2)),\]
\[b_2=\frac{1}{N}(\alpha_2^2 \beta_1 \gamma_1 +   \alpha_1^2 \beta_2 \gamma_2 -
\alpha_1 \alpha_2 (\beta_1 \gamma_1 +\beta_2 \gamma_2)
+ (\beta_1 - \beta_2)(\beta_1 \gamma_2 - \beta_2 \gamma_1) +\]\[
+ (\alpha_1 \alpha_2 (\alpha_1-\alpha_2)+(\alpha_1 \beta_2-\alpha_2 \beta_1))
(\delta_1 - \delta_2))  ,\]
\[b_3=\frac{1}{N}
((\alpha_2-\alpha_1) \beta_1 \beta_2 (\gamma_1 - \gamma_2) +
\alpha_1^2 \beta_2 \delta_1 + \alpha_2^2 \beta_1 \delta_2-
\alpha_1 \alpha_2 ( \beta_2 \delta_1 + \beta_1 \delta_2 )+\]
\[+(\beta_1 - \beta_2) (-\beta_2 \delta_1 + \beta_1 \delta_2)) \]
where $N$ is the resultant $(x_1-x_3)(x_1-x_4)(x_2-x_3)(x_2-x_4)$
or
\[  N = (\beta_1 -\beta_2 )^2 +(\alpha_1-\alpha_2)
(\alpha_1 \beta_2 -\alpha_2 \beta_1).\]
Because of (\ref{equation7}) we have
\[  A_3(x)=  x^2+(-\alpha_1-\alpha_2-\frac{a_0-2b_0 b_1}{b_0^2} )x
          +\frac{b_3^2-a_5}{b_0^2 \beta_1 \beta_2}=0         \]
and
\[  B_3(x)=-(y_5\frac{x-x_6}{x_5-x_6}+y_6\frac{x-x_5}{x_6-x_5}) =
-\frac{b(x_5)-b(x_6)}{x_5-x_6}x-\frac{b(x_6)x_5-b(x_5)x_6}{x_5-x_6}  \]
\[ =-(b_2 + b_1x_5 + b_1x_6 + b_0x_5^2 + b_0x_5x_6 + b_0x_6^2)x\]
\[    +b_3 + b_2x_5+ b_2x_6 +
  b_1x_5^2+ b_1x_6^2+ b_1x_5x_6 + b_0x_5^3+ b_0x_5^2x_6
+ b_0x_5x_6^2 + b_0x_6^3.  \]
Using $\alpha_3=-x_5-x_6$ and $\beta_3=x_5x_6$ we obtain
\[ B_3(x) =(-b_2+b_1 \alpha_3-b_0 \alpha_3^2+b_0 \beta_3)x
-b_0\alpha_3\beta_3+b_1\beta_3-b_3.  \]
Therefore we have the explicit rational group law
\begin{eqnarray*}
\alpha_3 & = & -\alpha_1-\alpha_2-\frac{a_0-2b_0 b_1}{b_0^2}, \\
 \beta_3 & = & \frac{b_3^2-a_5}{b_0^2 \beta_1 \beta_2},   \\
\gamma_3 & = & -b_2+b_1 \alpha_3-b_0 \alpha_3^2+b_0 \beta_3, \\
\delta_3 & = & -b_0\alpha_3\beta_3+b_1\beta_3-b_3  \\
\end{eqnarray*}
on the dense set of $Jac(C)-\Theta$ with
$(x_1-x_2)(x_3-x_4)Nb_0\beta_1\beta_2\neq 0$.

\noindent {\it Remark:}
The formulas are also true in the limit $x_1=x_2$, $x_3=x_4$.
The remaining special cases can be treated similar.


\begin{thebibliography}{9}


\bibitem[1]{Ca} Cantor, D.G., Computing in the Jacobian
                of a hyperelliptic curve,
                Mathematics of Computation {\bf 48},177(1987)95-101.

\bibitem[2]{Cas} Cassels, J.W.S.,
                The Mordell-Weil group of curves of genus 2.
                Arithmetic and geometry,
                Pap. dedic. I. R. Shafarevich,
                Vol. I: Arithmetic, Prog. Math. {\bf 35}(1983)27-60.

\bibitem[3]{Gr} Grant, D., Formal groups in genus two,
                J.reine angew. Math. {\bf 411}(1990)96-121.

\bibitem[4]{Hu} Huang, M.-D., Ierardi, D., Efficient Algorithms
                for the Riemann-Roch Problem and for Addition  in
                the Jacobian of a Curve,
                J. Symb. Comp. {\bf 18}(1994)519-539.

\bibitem[5]{J}  Jacobi, C.G.J.,
                \"Uber die Darstellung einer Reihe gegebener Werte
                durch eine gebrochene rationale
                Funktion, Crelle's J. {\bf 30}(1846)126-157.

\bibitem[6]{K}  Koblitz, N.,
                 Algebraic aspects of cryptography.
                 With an appendix on hyperelliptic curves,
                 Springer, New York 1999.

\bibitem[7]{M} Maseberg, S., Additionsformeln f\"ur
                Jacobi-Variet\"aten hyperelliptischer Kurven
                via Theta-Relationen, Diplomarbeit, Bremen 1998.

\bibitem[8]{Ma} Mazur, B., Arithmetic on curves,
                 Bull. AMS {\bf 14}(1986)207-259.

\bibitem[9]{Mu} Mumford, D.,
                 Tata lectures on theta.
                 Birkh\"{a}user, Boston 1994.

\bibitem[10]{St}  Staude, O.,
                 Geometrische Deutung der Additionstheoreme
                 der hyperelliptischen Integrale und Functionen erster Ordnung
                 im System der confocalen Fl\"{a}chen zweiten Grades,
                 Math. Ann. {\bf 22}(1883) 1-69, 145-176.

\bibitem[11]{Vol} Volcheck, E.J., Computing in the Jacobian of a Plane
                 Algebraic Curve,
                 in ANTS-I, Springer LNCS {\bf 877}(1994)221-233.

\bibitem[12]{Y} Yoshitomi, K., On height functions on Jacobian
                 surfaces, Manuscr. Math. {\bf 96}(1998)37-66.




\end{thebibliography}
\end{document}